\input AHTOH-E.STY

\UDC{%
515.142.331
+515.142.321
+512.543.16
%
%
}

\MSC{%
57K20, 
55U10,  
05E45,  
20F05
%
%
}

\title{%
%
Tiny non-Leighton complexes
}
\author{%
Natalia S. Dergacheva
\quad
Anton A. Klyachko
}
\address{%
Faculty of Mechanics and Mathematics of Moscow State University
\\
Moscow 119991, Leninskie gory, MSU.
\\
Moscow Center for Fundamental and Applied Mathematics
\\
nataliya.dergacheva@gmail.com
\quad
klyachko@mech.math.msu.su
}

\grants{\RSF 22-11-00075}

\abstract{%
How complex must two finite 2-complexes be to admit a common, but not 
finite common, covering?  We obtain an almost complete answer: the 
minimum possible number of triangles in a pseudo-simplicial triangulation 
of each complex is 3, 4, or 5.  }

\s 1.
Introduction

\proclaim Leighton's theorem {\rm[Lei82]}.
If two finite graphs have a common covering, then they have a common
finite covering.

\noindent
Alternative proofs and some generalisations of this result can be
found, e.g., in
[Neu10],
[BaK90],
[SGW19],
[Woo21],
[BrS21],
and references therein.

\noindent
For two-dimensional CW-complexes
(and cell coverings),
a similar theorem does not hold:
\-
the first example ([Wis96] and [Wis07]) of
such a \emph{non-Leighton}
pair of complexes
contains six two-cells
in each complex;
\-
later this number was reduced to four
[JaW09];
\-
and then to two [DK23].

\enditem
It is unknown
whether there exists a non-Leighton pair of complexes
with a single 2-cell (in both complexes,
or at least in one of them).
There exists a non-Leighton pair where one of the complexes
is \spacing{homeomorphic} to a complex with a single 2-cell
[DK24].

We consider another measure of complexity of CW-complexes:
the number of 2-simplices (triangles) in pseudo-simplicial
triangulations.
The \emph{pseudo-simplicial complexity} or simply \emph{complexity}
(see, e.~g., [JRT11]) of a 2-complex~$K$
is the minimal possible number of triangles
(2-cells) in a pseudo-simplicial complex homeomorphic to
$K$.
In other words, we consider
pseudo-simplicial triangulations of $K$
(this means that 2-cells
of a triangulation are triangles, but
two triangles are allowed to share several edges or vertices).
Thus, each $n$-gonal 2-cell ``costs"
$n-2$ (for $n\ge3$).

In this sense, the example from [JaW09] is better --- the
complexity (of both complexes) is 8;
while examples from [Wis07] and [DK23] are of complexity
12 and 24, respectively.

\proclaim Main theorem \rm(a simplified version).
There exist two finite 2-complexes of complexity five
with two 2-cells each
that have a common
covering, but have no finite common covering.

The explicit form of these complexes can be found at the end of Section~2.

Our approach is similar to that of [DK23].
The difference is that we succeed to use the Baumslag--Solitar
group~$\BS(1,2)$ instead of $\BS(3,5)$ used in [DK23];
it is kind of surprising:
\-
the arguments of [DK23] are based on
non-residual-finiteness of $\BS(3,5)$;
\-
replacing $\BS(3,5)$ with a smaller non-residually-finite group
$\BS(2,3)$ seems to be impossible, since
$\BS(2,3)$, unlike $\BS(3,5)$, contains the Klein-bottle group
as a subgroup,
which destroys the argument
(see the No-bottle lemma in [DK23]);
\-
in this paper, a trick allows us to replace
$\BS(3,5)$ with a residually finite (metabelian) group
$\BS(1,2)$.

\enditem
Nonetheless,
the fundamental groups of the both complexes comprising the non-Leighton
pair constructed in this paper are not residually finite.
The non-residual-finiteness plays a role in all known examples of
non-Leighton pairs.
In~[Wis07] and [JaW09],
the fundamental group of one of the two complexes comprising a
non-Leighton pair is residually finite (it is just the direct product of
two free groups).  However, the following question remains open.

\Question.
Does there exist a non-Leighton pair in which the fundamental groups of
both complexes are residually finite?

In Section 3, we prove the following fact.

\Proposition 1.
If two finite 2-complexes
have a common covering, but have no common finite covering,
then the complexity of each complex is
at least
three.

Thus, the minimum possible complexity of a member of a
non-Leighton pair is 3, 4, or 5.

{\noindent \bf Our notation}
is mainly standard. Note only that, if
$k\in \Z$, and $x$ and $y$ are elements of a group, then $x^y$, $x^{ky}$,
and $x^{-y}$ denote $y^{-1}xy$, $y^{-1}x^ky$, and $y^{-1}x^{-1}y$,
respectively;
The symbol $\Z_n$ (or $\gp x_n$) denotes the cyclic group of order $n$
(generated by $x$).
The symbol $F_2$ (or $F(x,y)$) denotes the free group of rank two
(with a basis $x,y$).
The symbol $*$ stands for the free product.
The \emph{Baumslag--Solitar groups} are
$\BS(n,m)\:=\pres<c,d|c^{nd}=c^m>$.
A semi-direct product is denoted $A\semitimes B$
(where the subgroup $B$ is normal).
The symbol $\epsilon$
denotes $\pm1$
throughout this paper.

\medskip

We thank Martin Bridson for a valuable remark
(see Section~3) and an anonymous referee for pointing out an error
in the original version of this paper.
The authors are grateful to
the Theoretical Physics and Mathematics Advancement Foundation ``BASIS".

\s 2.
Proof of the main theorem

Let us take
the standard presentations of the
fundamental groups of the torus and Klein bottle:
$$
G_1=\pres<a,b|ab=ba>\iso\BS(1,1)
\qqbox{and}
G_{-1}=\pres<a,b|ab=ba^{-1}>\iso\BS(1,-1),
$$
and consider the amalgamated free product
$H_\epsilon=G_\epsilon\*_{b=h}H$
(with a cyclic amalgamation) of
the groups $G_\epsilon$ and a group
$H=\pres<X|R>\supseteq\gp h_\infty$.
Let $K_\epsilon$ be the standard (one-vertex) complexes of
the following
presentations of $H_\epsilon$:
$$
H_\epsilon=
\pres<\{a\}\sqcup X|\left\{a^{\^h}=a^\epsilon\right\}\sqcup R>,
\qbox{where $\^h$ is a word in the alphabet $X^{\pm1}$
representing the element $h\in H$}.
$$
The Cayley graphs of the groups $G_\epsilon$ are surely isomorphic
(as abstract undirected graphs),
the same is true for the universal coverings of the standard
complexes of the presentations $G_\epsilon$
(these coverings are just the Euclidian plane partitioned into
squares).
A slightly less trivial observation is that,
for presentations $H_\epsilon$,
the universal coverings are isomorphic too:
\disp{\sl\narrower\narrower
for any group $H$ and any infinite-order element $h\in H$,
the universal coverings of the
complexes $K_\epsilon$ are isomorphic.
}
The readers may regard this fact as obvious,
but it
was explained in details in [DK23] (Observation~$(*)$).

If the complexes $K_\epsilon$ admit a common finite covering,
then the groups $H_\epsilon$ must
have isomorphic subgroups of \spacing{the same} finite index.
Indeed, a covering of degree $k$ of $K_\epsilon$
has $k$ vertices (therefore, isomorphic coverings must have the
same degree), and the fundamental group of the covering complex
is isomorphic to a subgroup of index $k$ in $H_\epsilon$.

Now, let us take a particular group
$
H=\BS(1,2)=\pres<c,d|c^d=c^2>
$
and a particular element $h=c\in H$.
To complete the proof of the main theorem,
it remains to prove the following fact.

\proclaim Weak-incommensurability lemma.
The groups $H_\epsilon=\pres<a,c,d|a^c=a^\epsilon,c^d=c^2>$,
where $\epsilon\in\{\pm1\}$,
contain no isomorphic subgroups of the
same finite index
\(i.e., if $U_{+1}\subseteq H_{+1}$ and $U_{-1}\subseteq H_{-1}$
are subgroups of the same finite index, then $U_{+1}\not\iso U_{-1}$\).

\Proof
Let $U\subseteq H_\epsilon$ be a finite-index subgroup.
Note that its intersection with the normal closure
$C\:=\nc c$ of the element $c\in H_\epsilon$
has the following properties:
\item{1)}
the index $|C:C\cap U|$ is odd;
\item{2)}
the subgroup $C\cap U$
is normal in $H_\epsilon$;
\item{3)}
$C\cap U=
\gp{\left\{u\in U\;\bigm|
\;\exists\;k\in\N\;\exists\;u'\in U\;u^{u'}=u^{2^k} \right\}}$.

\enditem
Indeed,
let $\^U\subseteq U$ be a normal in $H_\epsilon$
finite-index subgroup
(which always exists).
Let us choose $m$ such that $c^m\in\^U$
(for both $\epsilon=\pm1$; such $m$ exists because
$|H_\epsilon:\^U|<\infty$).
The integer $m$ can be assumed to be odd
(because $c$ and $c^2$ are conjugate).
Let $\phi$ be the
natural homomorphism
$$
H_\epsilon\too^\phi\~H_\epsilon\:=H_\epsilon/\nc{c^m}=
\cases{
(\gp d_\infty*\gp a_2)\semitimes\gp{\~c}_m&for $\epsilon=-1$
\cr\cr
F(d,a)\semitimes\gp{\~c}_m&for $\epsilon=+1$
}
\quad\pmatrix{\hbox{where the actions in the semidirect products}
\cr
\hbox{are the following: }
\~c^a=\~c,\;\~c^d=\~c^2}.
$$
\item{1)}
$\ker\phi=\nc{c^m}\subseteq C\cap\^U\subseteq C\cap U$,
hence,
$|C:C\cap U|=|\phi(C):\phi(C\cap U)|$,
where $\phi(C)=\gp{\~c}_m$ is a group of the odd order
$m$; therefore, the index is odd.
\item{2)}
$\ker\phi\subseteq C\cap U$ (see above), hence,
the normality of $C\cap U$ follows from the
normality of
$\phi(C\cap U)$, which is obvious:
$\phi(C)=\gp{\~c}_m$ is a normal cyclic subgroup,
therefore, all subgroups contained in it are normal
(all subgroups of a cyclic group are characteristic).
\item{3)}
{
The inclusion
$
C\cap U\supseteq
X\:=\gp{\left\{u\in U\;\bigm|
\;\exists\;k\in\N\;\exists\;u'\in U\;u^{u'}=u^{2^k} \right\}}
$
is easy to explain:
\itemitem{--}
$U\supseteq X$ by the definition of~$X$;
\itemitem{--}
the inclusion $C\supseteq X$
holds because $H_\epsilon/C$
has no nonidentity element conjugate to its even power,
see $(*)$ below).

Let us prove the reverse inclusion.
\itemitem{--}
$\gp c\cap U\subseteq X$.
Indeed, take $c^k\in\gp c\cap U$;
the index of~$U$ is finite,
hence, there is $q\in\Z\setminus\0$ such that $d^q\in U$;
so,
$(c^k)^{d^q}=(c^k)^{2^q}$ and $c^k\in X$
by the definition of~$X$;
\itemitem{--}
The same is true for subgroups conjugate to $\gp c$;
therefore,
$\ker\phi=\nc{c^m}\subseteq X$.
The inclusion $\nc{c^m}\subseteq C\cap U$ is explained above,
hence, $\ker\phi=\nc{c^m}\subseteq C\cap U\cap X$.
\itemitem{--}
Thus, it suffices to prove that
$\phi(C\cap U)\subseteq\phi(X)$.
This is true, because
$\phi(C)=\gp{\~c}$ and,
if $\phi(C)\ni\~c^k\in\phi(U)$,
then $c^k\in U\cdot\ker\phi=U$.
Hence, $c^k\in X$ (as $\gp c\cap U\subseteq X$,
see above).
Therefore, $\~c^k\in\phi(X)$, as required.

}

\noindent
Let us continue the proof of the lemma.
Let $U_{+1}$ and $U_{-1}$ be isomorphic subgroups of a finite index $i$
in~$H_{+1}$~and~$H_{-1}$.
Then
$\~U_\epsilon\:=U_\epsilon/(U_\epsilon\cap C)$
are isomorphic (by virtue of
Property 3)) subgroups of
$$
H_\epsilon/C=
\cases{
\gp d_\infty*\gp a_2&for $\epsilon=-1$
\cr
F(d,a)&for $\epsilon=+1$
}.
\eqno{(*)}
$$
the indices of these subgroups (of $\Z*\Z_2$ and $F_2$) are
$i/|C:(U_\epsilon\cap C)|$.
Indeed,
$$
i=|H_\epsilon:U_\epsilon|=
|H_\epsilon:U_\epsilon C|\cdot|U_\epsilon C:U_\epsilon|=
|H_\epsilon/C:\~U_\epsilon|\cdot|C:(U_\epsilon\cap C)|.
$$
Therefore,
the ratio
$|F_2:\~U_{+1}|/|\Z*\Z_2:\~U_{-1}|$
of these indices
has odd numerator and denominator
by virtue of Property 1).
But this is impossible, because
the subgroup $\~U_{+1}\subseteq F_2$ is free (by
the Nielsen--Schreier theorem)
and (by the Schreier formula)
$
\rk(\~U_{+1})-1=|F_2:\~U_{+1}|,
$
while, if $\~U_{-1}$ is free, then
$
\rk(\~U_{-1})-1={1\over2}|(\Z*\Z_2):\~U_{-1}|
$
%
by virtue of the following lemma.

\Lemma 1.
Let $U$ be a finite-index free subgroup of
$\Z*\Z_2$. Then \(its index is even and\/\)
$\rk(U)-1={1\over2}|(\Z*\Z_2):U|$.

\Proof
The kernel $F$ of the natural mapping $\Z*\Z_2\to\Z_2$
has index two and
rank two (it is freely generated by~$x$~and~$x^y$,
where $\Z*\Z_2=\gp x_\infty*\gp y_2$).
Then
$$
\eqalign{
&|(\Z*\Z_2):U|\cdot{\rk(U\cap F)-1\over\rk(U)-1}\=^S
\cr
&\=^S|(\Z*\Z_2):U|\cdot|U:(U\cap F)|=
|(\Z*\Z_2):(U\cap F)|=
|(\Z*\Z_2):F|\cdot|F:(U\cap F)|=
2\cdot|F:(U\cap F)|\=^S
\cr
&\=^S
2\cdot\bigl(\rk(U\cap F)-1\bigr)
\qbox{(where the equalities $\=^S$ are the Schreier formula).}
}
$$
Cancelling by $\rk(U\cap F)-1$
completes the proof of Lemma 1%
\fn{%
Lemma 1 is a special case of a general fact:
for each virtually free group, there is a well-defined
rational number called the \emph{virtual rank} [KZ24];
and these virtual ranks obey the Schreier formula.
}%
,
the weak-incommensurability lemma
(because the ratio of
indices $|F_2:\~U_{+1}|$ and $|(\Z*\Z_2):\~U_{-1}|$
must have odd numerator and denominator, as noted above),
and the following theorem.

\proclaim Main theorem.
The complexity-five complexes of the presentations
$$
H_\epsilon=\pres<a,c,d|
a^c=a^\epsilon,\ c^d=c^2>,
\qbox{where $\epsilon\in\{\pm1\}$,}
$$
containing two 2-cells
each
have a common covering, but have no finite common coverings.

\s 3.
Proof of Proposition 1

We call the minimal subcomplex of a 2-complex $K$
containing all its two-cells the \emph{two-dimensional
part} $K_2$ of $K$. Thus, $K=K_2\cup K_1$, where $K_1$
(the \emph{one-dimensional
part}) is a one-dimensional subcomplex
of $K$ such that $K_2\cap K_1$ is a 0-complex.

\proclaim Virtual-freedom lemma.
If a pseudo-simplicial 2-complex $P$ contains at most two 2-cells, then
either its fundamental group $\pi_1(P)$ is virtually free, or its
two-dimensional part $P_2$ contains a single vertex and is homeomorphic to
the torus or Klein bottle.

\Proof
Since contractions of edges joining two different vertices does
not change the fundamental group, we obtain that the fundamental group
of the complex has a presentation with at most two
relators of lengths at most three. Let us proceed with a simple
exhaustive search.

\-
If there is a letter occurring once in one of the relator
and not occurring in the other one, then
we apply the obvious Tietze transformation and obtain
a presentation with a single relator of length
at most three. Such a group is obviously virtually free
(it decomposes into a free product of a free group and a finite
group of order at most three).

\-
If there is a letter $a$ occurring several (two or three) times in one of
the relator and not occurring in the other one, then either we have
the case
considered above, or the presentation has a form
$$
G=\pres<a,b,c,\dots|a^{\pm2}b^{\pm1}=1,\;b^{\pm1}c^{\pm2}=1>
=\BS(1,-1)*F(d,e,\dots)
$$
or
$G=\pres<a,b,c,\dots|a^\alpha=1,\;b^\beta=1>$,
where $\alpha,\beta\in\{\pm2,\pm3\}$.
All the latter groups are virtually free.

\-
It remains to consider the case,
where each letter occurs in both
relators (or in none of them).
These are the following presentations
(up to obvious changes):
$\pres<a,b,c,\dots|a^\alpha=1,\;a^{\alpha'}=1>$,
$\pres<a,b,c,\dots|a^\alpha b^\beta=1,\;a^{\alpha'}b^{\beta'}=1>$
(where $\alpha,\alpha',\beta,\beta'\in\{\pm1,\pm2\}$)
and
$\pres<a,b,c,\dots|a^{\pm1}b^{\pm1}c^{\pm1}=1,\;a^{\pm1}b^{\pm1}c^{\pm1}=1>$.
All these groups are virtually free, except the last one,
where we have again the free product of a free group and
$\BS(1,\pm1)$.

\enditem
This completes the proof of the lemma.

\proclaim Bridson--Shepherd theorem
{\rm(cf. [BrS22], Corollary 3.6)}.
If two finite simplicial complexes have a common covering, and the
fundamental group of one of the complexes is virtually free, then the
complexes have a common finite covering.

In [BrS22], the condition is that the fundamental groups of \spacing{both}
complexes (instead of one of them) are \spacing{free} (instead of
virtually free). But the difference is inessential:

\-
if two complexes have a common covering, then their universal coverings
are isomorphic and, therefore, we can replace the initial complexes with
their finite coverings; but a complex with a virtually free fundamental
group has a finite covering with a free fundamental group;

\-
if the fundamental group of one of the complexes is (virtually) free, then
the common universal covering of these complexes is a quasi-tree;
therefore, the fundamental group of the other complex acts freely on this
quasi-tree; hence, this group is virtually free by the following fact
([Bu21], Corollary 5.8):
\disp{\sl
a finitely generated group acting freely on a locally finite graph which
is a quasi-tree is virtually free.
}%
As pointed out by Martin Bridson,
a shorter way of saying why the ``Bridson--Shepherd theorem'' above
is equivalent to its original version from [BrS22]
is the following:
the fundamental groups
of finite complexes
are quasi-isometric to their universal covers (more
precisely the 1-skeleta with the combinatorial metric), and being
virtually free is invariant with respect to quasi-isometry.

\proclaim Gluing lemma.
Let $K'$ and $K''$ be subcomplexes of a finite complex
$K=K'\cup K''$, and let $L'$ and $L''$
be subcomplexes of a finite complex
$L=L'\cup L''$. Suppose that
\item{\rm1)}
$K'\cap K''$ is a single vertex $v$,
and $L'\cap L''$ is
a 0-complex consisting of
vertices $w_1,\dots,w_q$;
\item{\rm2)}
$K'$ and $L'$ have a common finite
covering $K'\oot^{\kappa'}M'\too^{\lambda'}L'$,
consistent on the
intersection:
$$
(\kappa')^{-1}(K'\cap K'')=(\lambda')^{-1}(L'\cap L'');
$$
\item{\rm3)}
$K''$ and $L''$ have a common
finite covering $K''\oot^{\kappa''}M''\too^{\lambda''}L''$,
consistent on the intersection:
$$
(\kappa'')^{-1}(K'\cap K'')=(\lambda'')^{-1}(L'\cap L'').
$$
\enditem
Then
$K$ and $L$ have a common finite covering.

\Proof
Suppose that the covering $\kappa'$ has degree $k'$
and the covering
$\kappa''$ has degree $k''$.
Then, the degrees of the coverings $\lambda'$
and $\lambda''$
are $k'/q$ and $k''/q$, respectively,
because of the consistency.

Take the disjoint union
$M'_1\sqcup\dots\sqcup M'_{k''}\sqcup M''_1\sqcup\dots\sqcup M''_{k'}$
(where $M'_i$ and $M''_j$ are
copies of complexes $M'$ and $M''$).
For each $j\in\{1,\dots,q\}$,
we identify in these complexes $k'k''/q$ preimages
of the vertex $w_j$ lying in
$M_i'$ with~$k'k''/q$ preimages of the vertex $w_j$ lying in $M_i''$.
Take the connected component $X$ of the obtained finite complex and
the natural covering $K\ot X\to L$. This completes the proof of the lemma.


\proclaim Leighton's theorem for coloured graphs \rm[Neu10]
(a special case%
\fn{\rm
The general case consists in that we can colour
not only vertices, but also edges.
[Neu10] contains a proof and also a note
that,
first, this coloured theorem easily follows from the classical
Leighton theorem, and secondly, the original argument of Leighton [Lei82]
proves, in fact, the coloured theorem.}%
).
If two finite graph $A$ and $B$ whose vertices are coloured with some
colours, have a common covering $C$ and
the mappings $A\ot C\to B$
preserve colours,
then $A$ and $B$ have a common finite covering, and the corresponding
mappings preserve colours too.

\medskip

Let us proceed with the proof of Proposition 1.
Suppose that finite pseudo-simplicial two-complexes $K$ and $L$
have a common covering, and
$K$ has at most two
2-cells.
We want to find a finite common covering.

The Bridson--Shepherd theorem and
the virtual-freedom lemma reduce the situation to
the case, where
the two-dimensional part~$K_2$ of $K$
is homeomorphic to the torus or Klein bottle and
contains a single vertex $v$.

The common covering $K\ot M\to L$ of $K$ and $L$ gives a common
covering $K_2\ot M_2\to L_2$ of $K_2$~and~$L_2$, and also a common
covering $K_1\ot M_1\to L_1$ of $K_1$~and~$L_1$. This means that
\item{2)}
$K_2$ and $L_2$ have a common finite covering
$K_2\ot P\to L_2$,
and $L_2$ is also
homeomorphic to the torus or the Klein bottle
(and contains vertices $w_1,\dots,w_q$);
\item{1)}
for complexes $K_1$ and $L_1$, there is a common finite covering
$K_1\ot R\to L_1$, and the full preimage of $v$
coincides with the full preimage of the set $\{w_1,\dots,w_q\}$.

\enditem
Indeed,
\item{2)}
the universal covering of the torus or Klein bottle
(i.e., the plane) must cover
$L_2$, i.e. $L_2$ is  homeomorphic to a surface,
this surface must be the torus or Klein bottle 
(because the fundamental groups of other surfaces
are not commensurable with $\Z\oplus\Z$);
\item{1)}
the preimages of $v$ and $w_i$ in $M$ (and in $M_1$) must coincide
(because this is precisely $M_2\cap M_1$); therefore, we can
colour $v$, $w_i$, and all their preimages in $M_1$ red,
colour all the remaining vertices of $K_1$, $L_1$, and $M_1$ black, and
apply Leighton's theorem for coloured graphs.

\enditem
To complete the proof, it remains to apply the gluing lemma
(putting $K'=K_1$, $L'=L_1$, $K''=K_2$, and $L''=L_2$).


\References

[BaK90]
H. Bass, R. Kulkarni,
Uniform tree lattices,
J. Amer. Math. Soc., 3:4 (1990), 843-902.

[BrS22]
M. Bridson, S. Shepherd,
Leighton's theorem: extensions, limitations, and quasitrees,
Algebraic and Geometric Topology, 22:2 (2022), 881-917.
\arXiv 2009.04305.


[Bu21]
J. O. Button,
Groups acting on hyperbolic spaces with a locally finite orbit.
arXiv: 2111.13427 (2021).

[DK23]
N. S. Dergacheva, A. A. Klyachko,
Small non-Leighton two-complexes,
Math. Proc. Cambridge Philos. Soc., 174:2 (2023), 385-391.
\arXiv 2108.01398

[DK24]
N. S. Dergacheva, A. A. Klyachko,
Forester's lattices and small non-Leighton complexes,
Journal of Topology and Analysis (to appear),
\doi 10.1142/S1793525324500547
\arXiv:2407.06680

[KZ24]
A. A. Klyachko, A. O. Zakharov,
An analogue of the strengthened Hanna Neumann conjecture
for virtually free groups and virtually free products,
Michigan Mathematical Journal (in publication).
\arXiv 2106.05821

[JaW09]
D. Janzen, D. T. Wise,
A smallest irreducible lattice in the product of trees,
Algebraic and Geometric Topology, 9:4 (2009), 2191-2201.

[JRT11]
W. Jaco, J. H, Rubinstein, S. Tillmann,
Coverings and minimal triangulations of 3-manifolds,
Algebraic \& Geometric Topology, 11:3 (2011), 1257-1265.
\arXiv:0903.0112

[Lei82]
F. T. Leighton,
Finite common coverings of graphs,
J. Combin. Theory, Series B, 33:3 (1982), 231-238.

[Neu10]
W. D. Neumann,
On Leighton's graph covering theorem,
Groups, Geometry, and Dynamics, 4:4 (2010), 863-872.
\arXiv 0906.2496

[SGW19]
S. Shepherd, G. Gardam,  D. J. Woodhouse,
Two generalisations of Leighton's Theorem,
arXiv:1908.00830.

[Wis96]
D. T. Wise,
Non-positively curved squared complexes:
Aperiodic tilings and non-residually finite groups.
PhD Thesis, Princeton University, 1996.

[Wis07]
D. T. Wise,
Complete square complexes,
Commentarii Mathematici Helvetici, 82:4 (2007), 683-724.

[Woo21]
D. Woodhouse,
Revisiting Leighton's theorem with the Haar measure,
Mathematical Proceedings of the Cambridge Philosophical Society,
170:3 (2021), 615-623.
\arXiv 1806.08196

\end